\documentclass{amsart}

\usepackage{graphicx}

\usepackage{amsmath}

\usepackage{amsfonts}

\usepackage{amssymb}

\newtheorem{theorem}{Theorem} [section]

\newtheorem{thm}[theorem]{Theorem}

\newtheorem{prop}[theorem]{Proposition}

\begin{document}

\title{Intrinsically Linked Graphs with Knotted Components}

\author{Thomas Fleming}

\begin{abstract}

We construct a graph $G$ such that any embedding of $G$ into $R^{3}$
contains a nonsplit link of two components, where at least one of
the components is a nontrivial knot.  Further, for any $m < n$ we
produce a graph $H$ so that every embedding of $H$ contains a
nonsplit $n$ component link, where at least $m$ of the components
are nontrivial knots.

We then turn our attention to complete graphs and show that for any
given $n$, every embedding of a large enough complete graph contains
a two component link whose linking number is a nonzero multiple of
$n$.

\end{abstract}

\maketitle

\section{Introduction}

We call a graph \emph{intrinsically linked} (\emph{intrinsically
knotted}) if every embedding of that graph into three-space contains
a nontrivial link (knot). Intrinsically linked and intrinsically
knotted graphs were first studied in the early 1980s by Sachs
\cite{sachs} and by Conway and Gordon \cite{c&g}.  Since then these
graphs have been extensively studied. Significant progress has been
made in this area, such as the complete classification of minor
minimal examples for intrinsically linked graphs by Robertson,
Seymour, and Thomas \cite{r&s}.  Since this classification was
completed, work has turned to finding graphs whose every embedding
contains more complex structures.  For example, one could require
that every embedding of the graph contain a nonsplit link of $n$
components \cite{flapan2}, or a link with linking number larger than
some constant \cite{flapan3, s-t}.

We would now like to combine the two properties of intrinsic linking
and intrinsic knotting and ask the following question:  Does there
exist a graph such that any embedding of this graph contains a
nonsplit link, where at least one of the components of this link is
a nontrivial knot?

The answer is affirmative, and we demonstrate examples.

\medskip

\noindent \textbf{Theorem \ref{mainthm}} \emph{Every embedding of
the graph $F(126)$ contains a nonsplit link of two components, where
at least one of these two components is a nontrivial knot.}

\medskip

This graph contains one hundred and twenty-six vertices, so it is
unlikely to be minor minimal. In fact, it is possible that graphs as
small as $K_{10}$ or $K_{11}$ have this property, but current
techniques do not allow us to detect it.

A graph containing an $n$-component link in every embedding was
first demonstrated by Flapan, Foisy, Naimi, and Pommersheim
\cite{flapan2}. The link constructed in \cite{flapan2} has the form
of a chain of components, each linked to its neighbors.  A different
type of $n$-component link, more like a ring of keys where one
component links all the others, was found in \cite{mend}.  Using
this later construction, we can produce a link
 of arbitrarily many components, some of which are nontrivial knots.
\medskip

\noindent \textbf{Theorem \ref{knottedkeys}} \emph{Given $m \leq n$,
there exists a graph $G$ such that every embedding of $G$ contains a
nonsplit link of $n+1$ components, where at least $m$ of the
components are nontrivial knots.}

\medskip

Flapan \cite{flapan3} and Shirai and Taniyama \cite{s-t} proved that
for a given $k$, certain graphs always contain a two component link
with linking number greater than or equal to $k$. This led the
author and Diesl to study graphs that contained two component links
with linking number equal to a multiple of $k$, for $k =2^{r}$
\cite{mend}. Here we continue this work and search for links with
linking number equal to a multiple of a given integer $n$.

Let $\eta_{n} = \alpha'_{n}(\zeta_{n}+3)$ where $\zeta_{n} =
(n+1)(\lfloor \frac{n+1}{2} \rfloor)(\lceil \frac{n+1}{2}
\rceil)^{\lceil log_{2}(n) \rceil-2}$, and $\alpha'_{n}$ is the
sequence from \cite{mend} defined as $\alpha_{1}' =1$, $\alpha_{2m}'
= 2\alpha_{2m-1}'$, and $\alpha_{2m+1}'=2\alpha_{2m}'+1$. Then we
have the following theorem.
\medskip

\noindent \textbf{Theorem \ref{zeromodn}} \emph{Given $n \geq 5$,
every embedding of $K_{\eta_{n}}$ contains a 2-component link $L,Z$
with $lk (L,Z) = kn$ for some $k \neq 0$.}

\medskip

The author and Diesl produced the analogous result for the case of
$n=2^{r}$ and $n=3$ in \cite{mend}.  When $n=2^{r}$ the number of
vertices required for the graph in Theorem \ref{zeromodn} is roughly
$\alpha'_{n}2^{r^{2}}$ whereas the bound from \cite{mend} is
smaller, roughly $\alpha_{n}'2^{\frac{r^{2}}{2}}$. So, while here we
are able to extend the construction to all $n$, this extension comes
at a cost.

\section{Knotted Links}

To prove Thoerem \ref{mainthm}, we first produce an intrinsically
knotted graph $F$ and show that a certain pair of edges in this
graph are contained in a knotted cycle in any embedding.  This graph
$F$ is closely related to the graph introduced by Foisy in
\cite{foisy2}.

We then take ten copies of $F$ and carefully glue them together to
produce $F(126)$.  The crucial idea is that any embedding of
$F(126)$ can be contracted down to a copy of $K_{6}$, where every
triangle containing a distinguished vertex is the contraction of a
knotted cycle.  Thus the linked triangles in $K_{6}$ correspond to a
pair of linked cycles in $F(126)$ and as one of these triangles
contains the distinguished vertex, the corresponding cycle is
knotted.

\begin{thm}Every embedding of the graph $F(126)$ contains a nonsplit link of two components, where at least one of these two components is a nontrivial knot.

\label{mainthm}

\end{thm}

\begin{proof}

Let $F$ be the graph shown in Figure \ref{fgraph}.  Note that this
graph contains the graph produced by Foisy in \cite{foisy2} as a
minor, and hence is intrinsically knotted.

\begin{figure}[hbtp]

\centering

\includegraphics[width=.4\textwidth]{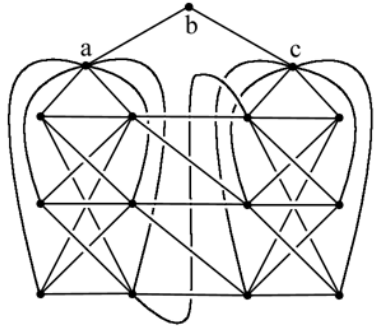}

\caption{The graph $F$.}

\label{fgraph}

\end{figure}

We, however, want the edges $a-b$ and $b-c$ to always be contained
in a knotted cycle, whatever the embedding of $F$. This graph
contains two disjoint copies of $K_{3,3,1}$ which is intrinsically
linked, and in fact each copy always contains a triangle linked with
a square. Using the edges running between the two copies of
$K_{3,3,1}$, we may contract $F$ to the graph shown in Figure
\ref{4cycles}.  We contract the triangle containing $a$ to cycle 1,
and the triangle containing $c$ to cycle 2.  Because of the edges in
the middle of the graph, the triangle containing $a$ is adjacent to
the square linked with the triangle containing $c$. Contract this
square to cycle 2, and the analogous square in the other copy of
$K_{3,3,1}$ to cycle 3.  This gives us that cycle 1 is linked with
cycle 3, and cycle 4 is linked with cycle 2.  It is shown in
\cite{foisy} and \cite{ty} that the graph in Figure \ref{4cycles}
always contains a knotted cycle, and clearly this knotted cycle must
use the edges $a-b$ and $b-c$. Since vertex expansion (and edge
addition) do not change the isotopy type of a cycle, every embedding
of $F$ contains a knotted cycle that uses the edges $a-b$ and $b-c$.
(Note that this would not be true if we used the corresponding
expansion of Foisy's graph from \cite{foisy2}).

\begin{figure}[hbtp]

\centering

\includegraphics[width=.4\textwidth]{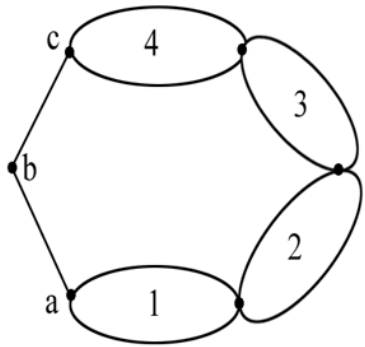}

\caption{A contraction of $F$. The cycles labeled 1 and 4 come from the triangles containing a and
c, respectively. Cycles 2 and 3 are the squares linked with the triangles corresponding to 4 and 1
respectively.} \label{4cycles}

\end{figure}

We now take ten copies of $F$.  Label the upper vertices of the first graph $a_{1}, b_{1}, c_{1}$ and the upper vertices of the second graph $a_{2}, b_{2}, c_{2}$ and so on.

Take the first four copies and identify all the edges $a_{i}-b_{i}$.  Note that we have six labeled vertices after this identification.  The resulting graph appears in Figure \ref{glue1}.  The circles containing the letter $F$ denote the undrawn remainder of the graph $F$.

\begin{figure}[hbtp]

\centering

\includegraphics[width=.4\textwidth]{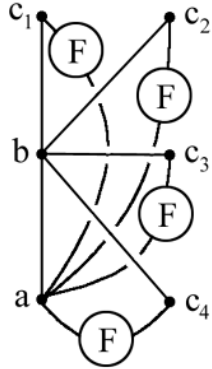}

\caption{Four copies of $F$ glued together. The edge labeled $a-b$ is the result of identifying the $a_{i}-b_{i}$.}

\label{glue1}

\end{figure}

Take the remaining six copies of $F$ and glue them on as follows.  Identify $a_{5}-b_{5}-c_{5}$ to $c_{1}-b-c_{2}$. Identify $a_{6}-b_{6}-c_{6}$ to $c_{1}-b-c_{3}$. Identify $a_{7}-b_{7}-c_{7}$ to $c_{1}-b-c_{4}$. Identify $a_{8}-b_{8}-c_{8}$ to $c_{2}-b-c_{3}$. Identify $a_{9}-b_{9}-c_{9}$ to $c_{2}-b-c_{4}$. Finally, identify $a_{10}-b_{10}-c_{10}$ to $c_{3}-b-c_{4}$.

This graph is $F(126)$.

Choose an embedding of $F(126)$.  In each copy of $F$, the edges
$a-b$ and $b-c$ are used in a knotted cycle.  Contract each copy of
$F$ along this knotted cycle, deleting any edges that become
parallel to it.  The result is shown in Figure \ref{k6}.  Note that
we now have a copy of $K_{6}$ such that every triangle that includes
vertex $b$ is the contraction of knotted cycle, and hence knotted.

\begin{figure}[hbtp]

\centering

\includegraphics[width=.4\textwidth]{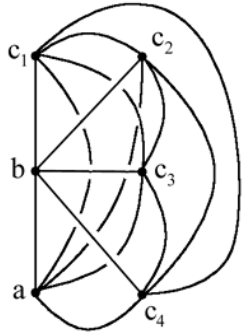}

\caption{The contraction to $K_{6}$. Note that any triangle that uses two straight lines and one curved one is knotted.}

\label{k6}

\end{figure}

Every embedding of $K_{6}$ contains a pair of linked triangles, and one of these triangles must include the vertex $b$.  Thus in this embedding of $K_{6}$ we have a two component link, where the component using vertex $b$ is knotted. Again, as vertex expansion and edge addition do not change the isotopy class of these cycles, we have that every embedding of $F(126)$ contains two linked cycles, where at least one of them is a nontrivial knot.  This completes the proof.

\end{proof}

Using the techniques of Theorem \ref{mainthm}, we can produce a slightly smaller graph with the same property as $F(126)$. We will follow the proof of Theorem \ref{mainthm} exactly, except we now glue nine copies of $F$ together to produce a graph that contracts to $K_{3,3,1}$.  A similar argument may be possible for other graphs in the Peterson family.

\begin{prop} Every embedding of $F(115)$ contains a nonsplit two component link, where at least one of these components is a nontrivial knot.

\end{prop}

\begin{proof}

We will glue nine copies of $F$ together to produce a graph with one
hundred and fifteen vertices that contracts down to $K_{3,3,1}$.

To construct $F(115)$, take nine copies of $F$, labeling their
vertices as before.  Identify $b_{1}, b_{2}, b_{3}$ and call the
resulting vertex $B$.  We now have seven preferred vertices, $a_{1},
a_{2}, a_{3}, B, c_{1}, c_{2}, c_{3}$.  These will form the
$K_{3,3,1}$.

Now identify $a_{4}-b_{4}-c_{4}$ to $a_{1}-B-c_{2}$. Identify $a_{5}-b_{5}-c_{5}$ to $a_{1}-B-c_{3}$. Identify $a_{6}-b_{6}-c_{6}$ to $a_{2}-B-c_{1}$. Identify $a_{7}-b_{7}-c_{7}$ to $a_{2}-B-c_{3}$. Identify $a_{8}-b_{8}-c_{8}$ to $a_{3}-B-c_{1}$. Identify $a_{9}-b_{9}-c_{9}$ to $a_{3}-B-c_{2}$.

We proceed as above, contracting along knotted cycles to obtain a copy of $K_{3,3,1}$ where every triangle containing $B$ is a knotted cycle.  Every embedding of $K_{3,3,1}$ contains a triangle linked with a square, and the triangle must contain $B$.  This gives the desired result.

\end{proof}

%\begin{cor}There exists a graph such that every embedding contains a nonsplit %three component link, where at least two of the components are nontrivial %knots.

%\end{cor}

%\begin{proof}

%To construct such a graph, we start with two copies of $F(126)$.  The %potentially unknotted component of the link in $F(126)$ contains three of $a, %c_{1}, c_{2}, c_{3}, c_{4}$ and hence at least two of  the $c_{i}$. Add an edge %from $c_{1}$ in the first copy to $c_{1}$ in the second copy.  Repeat for the %other $c_{i}$.

%Now, as the two linked components of $K_{6}$ have nonzero linking number, the %two linked cycles produced in the proof of Theorem \ref{mainthm} have nonzero %linking number in every embedding.  We may now apply the gluing Lemma 3 of %\cite{foisy3} to the potentially unknotted component in each copy of $F(126)$.  %This produces a cycle that has nonzero linking number with the knotted %component in each copy of $F(126)$.

%\end{proof}

By combining many copies of these graphs, we can construct a graph whose every embedding contains a nonsplit link with arbitrarily many knotted components. In fact, our construction produces an $n+1$-component link where all but one of the components are knots.

\begin{thm}
Given $m \leq n$, there exists a graph $G$ such that every embedding
of $G$ contains a nonsplit link of $n+1$ components, where at least
$m$ of the components are nontrivial knots.

\label{knottedkeys}
\end{thm}

\begin{proof}

By Theorem \ref{mainthm}, every embedding of the graph $F(126)$ contains a two component link $L, K$ where $K$ is a nontrivial knot and $lk(L,K) = 1$ mod $2$.

We may now apply Corollary 2.4 of \cite{mend}, which states that if
every embedding of $H$ contains a two component link $L,Z$ with
non-zero linking number, then every embedding of $*^{\alpha_{n}'}H$
contains an $n+1$ component link $L, Z_{i}$ with $lk(L,Z_{i}) \neq
0$. Here $\alpha_{1}' =1$, $\alpha_{2m}' = 2\alpha_{2m-1}'$, and
$\alpha_{2m+1}'=2\alpha_{2m}'+1$.

%Note that $\alpha_{2m-1}' = \frac{4^{m}-1}{3}$.

In the proof of this corollary, the desired link in an embedding of
$*^{\alpha_{n}'}H$ is produced by inductively modifying cycles in
pairs of links with $k-1$ components to produce a link with $k$
components.  However, it is only necessary to modify the central
component $L$ at each stage, so the $Z_{i}$ remain unchanged
throughout the process.

Choose $H = F(126)$.  In $G=*^{\alpha_{n}'}F(126)$, there are
$\alpha_{n}'$ copies of $F(126)$. In the $i$th copy of $F(126)$, we
can find a link $L_{i}, K_{i}$ by Theorem \ref{mainthm}.  Thus, in
$G=*^{\alpha_{n}'}F(126)$, we can find a nonsplit $n+1$ component
link $L, Z_{i}$, where the $Z_{i}$ are the cycles $K_{i}$ and the
cycle $L$ is obtained from the $L_{i}$.  The $K_{i}$ are nontrivial
knots, so we are done.

\end{proof}

\section{Linking modulo n}

In \cite{mend} the author and Diesl were able to show that for large
enough $m$, every embedding of $K_{m}$ contains a two component link
with linking number $k2^{r} (k \neq 0)$.  Here we extend this
construction from powers of two to all $n$.

Crucial to these constructions is the ``ring of keys'' lemma (Lemma
2.2 of \cite{mend}).  In any embedding of a large enough complete
graph, this lemma ensures that we can find an $n+1$ component link
$L, Z_{i}$ where the preferred component $L$ has non-zero linking
number with the other $n$ components. Further, as discussed in the
proof of Theorem 2.5 of \cite{mend}, using a result of Johnson and
Johnson \cite{johnson}, this link may be constructed so that the
number of vertices in each of the $Z_{i}$ is bounded from below.

Let $\alpha'_{2m-1} = \frac{4^{m} -1}{3}$, and $\alpha'_{2m} =
2\alpha'_{2m-1}$, then Corollary 2.4 of \cite{mend} implies that if
$G$ is intrinsically linked, we may use Lemma 2.2 \cite{mend} to
show every embedding of $*^{\alpha'_{n+1}}G$ contains an $n+1$
component link of the desired type.  Some of the properties of the
components $Z_{i}$ are inherited from $G$; it is this fact that is
useful in the proof of Theorem \ref{knottedkeys}.

For  $n \geq 5$, let $\eta_{n} = \alpha'_{n}(\zeta_{n}+3)$ where
$\zeta_{n} = (n+1)(\lfloor \frac{n+1}{2} \rfloor)(\lceil
\frac{n+1}{2} \rceil)^{\lceil log_{2}(n) \rceil-2}$.  This
$\zeta_{n}$ is the bound on the minimal number of vertices in the
link components $Z_{i}$.

For $n = 2,3,4$, the theorem below is true and the argument proceeds
with fewer vertices required in $K_{\eta_{n}}$. However, the bounds
produced by this result for $n =2,3,4$ are inferior to those of
\cite{mend}.

\begin{thm}Given $n \geq 5$, every embedding of $K_{\eta_{n}}$ contains a 2-component link $L,Z$ with
$lk (L,Z) = kn$ for some $k \neq 0$. \label{zeromodn}
\end{thm}

\begin{proof}

Given $n$, by Lemma 2.2 of \cite{mend}, in any embedding of
$K_{\eta_{n}}$ there exist disjoint cycles $L, Z_{1}, Z_{2} \ldots
Z_{n}$ with $lk(L,Z_{i}) \neq 0$ and each $Z_{i}$ containing at
least $\zeta_{n}$ vertices. Choose the orientations of $Z_{i}$ so
that this linking number is positive. Examine the sums $S_{j} :=
\sum_{i=1}^{j} [Z_{i}]$ in $H_{1}(\mathbb{R}^{3} \setminus L)$.  As
there are $n$ such sums, either two are equal modulo $n$, or one is
zero modulo $n$.  In the latter case this sum is $kn$ for $k \neq 0$
since the $[Z_{i}]$ are all positive.  In the former case, say
$S_{j_{1}} \equiv S_{j_{2}}$ and $j_{1} < j_{2}$.
$$\sum_{j_{1}+1}^{j_{2}} [Z_{i}] \equiv \sum_{1}^{j_{1}} [Z_{i}] +
\sum_{j_{2}+1}^{n} [Z_{i}] \equiv 0 \quad \mbox{mod n}$$ Thus, there
must be some set $J$ of the $Z_{i}$ with $\sum_{i \in J} lk(L,Z_{i})
= kn$ with $k \neq 0$.
%We will induct on the size of $J$.

If $J$ is size one, then we are done.  We now induct on the size of
$J$.

Suppose $J$ contains $r$ cycles.  Select two of them, say $Z_{1}$
and $Z_{2}$, and select $n+1$ vertices evenly spaced around $Z_{i}$.
Since $K_{\eta_{n}}$ is complete, each preferred vertex is adjacent
to all the others. Label the vertices of $Z_{1}$ cyclically in an
order agreeing with the orientation of $Z_{1}$ and label the
vertices of $Z_{2}$ cyclically counter the orientation of $Z_{2}$.
We form cycles $A_{i}$ by beginning at vertex $i$ on $Z_{1}$, taking
the path along $Z_{1}$ to vertex $i+1$, then the edge from vertex
$i+1$ on $Z_{1}$ to vertex $i+1$ on $Z_{2}$, the path along $Z_{2}$
to vertex $i$, and then the edge from vertex $i$ on $Z_{2}$ to
vertex $i$ on $Z_{1}$. Note that by construction we have $n+1$ such
cycles and $\sum [A_{i}] = [Z_{1}] + [Z_{2}] \in
H_{1}(\mathbb{R}^{3} \setminus L)$. See Figure \ref{aloops}.

\begin{figure}[hbtp]
\centering
\begin{picture}(324,240)
\includegraphics{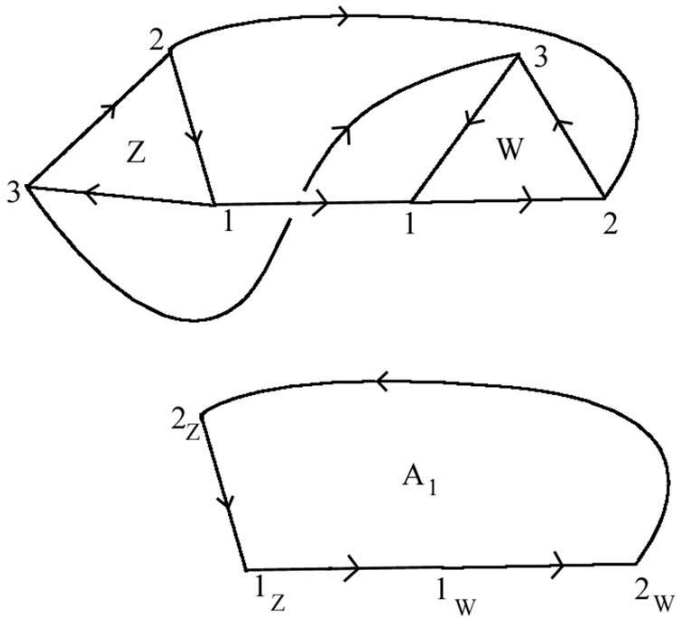}
\end{picture}
\caption{Constructing the $A_{i}$ for n=2}
\label{aloops}
\end{figure}

%Need to form $n+1$ cycles, not $n$.  This will inflate the count a
%bit!

 Using the same logic as for selecting the set $J$ we now
examine the $A_{i}$. Form the partial sums $T_{j} := \sum_{i=1}^{j}
[A_{i}]$ for $1 \leq j \leq n$. If $T_{j} \equiv [Z_{1}] + [Z_{2}]
\quad \mbox{mod n}$, then we remove $Z_{1}$ and $Z_{2}$ from $J$ and
replace them with $\cup_{1}^{j} A_{i}$.  We have reduced the size of
$J$ by one, but have not altered the sum of the elements in
$H_{1}(\mathbb{R}^{3} \setminus L)$.

If no partial sum equals $[Z_{1}] + [Z_{2}]$ modulo $n$, then
$T_{j_{1}} \equiv T_{j_{2}}$ for some $j_{1} < j_{2}$.  Then
$\sum_{j_{1}+1}^{j_{2}} [A_{i}] \equiv  0$ modulo $n$ so
$\sum_{1}^{j_{1}} [A_{i}] + \sum_{j_{2}+1}^{n+1} [A_{i}] \equiv
[Z_{1}] + [Z_{2}]$ modulo $n$. Let $A =
\cup_{j_{1}+1}^{j_{2}}[A_{i}]$ and $Z_{1}' = \cup_{1}^{j_{1}}
[A_{i}] \cup \cup_{j_{2}+1}^{n+1} [A_{i}]$. Note that $Z_{1}'$ is a
single cycle, as $A_{1}$ and $A_{n+1}$ are adjacent.  If $lk(A,L)
\neq 0$ then we are done. If $lk(A,L) = 0$ then $lk(Z'_{1},L) =
lk(Z_{1},L) + lk(Z_{2},L)$. So, once again, we may remove $Z_{1}$
and $Z_{2}$ from $J$ and replace them with $Z_{1}'$.

Since we may choose which cycles to pair at each step, we may pair
the cycles in $J$ with the largest number of vertices. Note that if
$Z_{1}$ and $Z_{2}$ have $(n+1)m$ vertices, the $A_{i}$ will contain
at least $2m+2$ vertices and hence the cycle $Z_{1}'$ will contain
at least $2m+2$ vertices. In the worst case, we must continue this
pairing process until only a single cycle remains in $J$. Thus a
cycle and its descendants will be used in at most $\lceil log_{2}(n)
\rceil$ pairings.  The initial cycles $Z_{i}$ contain at least
$(n+1)(\lfloor \frac{n+1}{2} \rfloor)(\lceil \frac{n+1}{2}
\rceil)^{\lceil log_{2}(n) \rceil -2}$ vertices, so after pairing
the elements in $J$ and producing the new cycles, $J$ will have
$\lceil \frac{n}{2} \rceil$ elements each with at least $2(\lfloor
\frac{n+1}{2} \rfloor)(\lceil \frac{n+1}{2} \rceil)^{\lceil
log_{2}(n) \rceil -2} + 2 \geq (n+1)(\lfloor \frac{n+1}{2}
\rfloor)(\lceil \frac{n+1}{2} \rceil)^{\lceil log_{2}(n) \rceil -3}
\geq (n+1)(\lfloor \frac{n+1}{2} \rfloor)(\lceil \frac{n+1}{2}
\rceil)^{\lceil log_{2}(\frac{n}{2}) \rceil -2}$ vertices.  The set
$J$ is now smaller, and the elements of $J$ have sufficiently many
vertices, so we may apply the induction hypothesis.

%In the case $n=2$, three vertices are required (cite flapan?), and
%when $n =3,4$ only $n(\lfloor \frac{n}{2} \rfloor)$ vertices are
%required.

\end{proof}

%Remark -- compare this to the bounds in AGT paper ($\beta_{r}$). The
%size of alpha required is the same, but here I will have roughly
%$2^{r^{2}}$ +3 whereas the the diesl paper, gamma is about
%$2^{r^{2}/2}$ (triangle number in exponent. Thus, mend bound is
%better for $2^r$, but that's good. can make a case for both papers
%then. This one does all n, but that one does powers of two more
%efficiently.

%Do I have the correct value for $K_{tba}$?  need to go through and
%correct any count that has $n$ in it!

\medskip

\textsc{University of California San Diego, Department of Mathematics, 9500 Gilman Dr., La Jolla,
CA 92093-0112}

\emph{E-mail address:} \texttt{tfleming@math.ucsd.edu}
\end{document}